\documentclass[12pt,american]{article}
\usepackage[T1]{fontenc}
\usepackage[latin9]{inputenc}
\usepackage{color}
\usepackage{babel}
\usepackage{varioref}
\usepackage{mathrsfs}
\usepackage{amsmath}
\usepackage{amsthm}
\usepackage{amssymb}
\usepackage{graphicx}
\usepackage[unicode=true,
 bookmarks=true,bookmarksnumbered=true,bookmarksopen=false,
 breaklinks=true,pdfborder={0 0 1},backref=false,colorlinks=true]
 {hyperref}
\hypersetup{pdftitle={Asymptotic Behavior of the Subordinated Traveling Waves},
 pdfauthor={{Jos{\'e} L.\ da Silva and Yuri G. Kondratiev},
 pdfsubject={Asymptotic Behavior of the Subordinated Traveling Waves},
 pdfkeywords={General fractional derivative, Convolution operator, Subordination principle},
 },linkcolor=black,citecolor=black,urlcolor=black,filecolor=black}

\makeatletter
\theoremstyle{plain}
\newtheorem{thm}{\protect\theoremname}
\theoremstyle{definition}
\newtheorem{example}[thm]{\protect\examplename}
\theoremstyle{plain}
\newtheorem{prop}[thm]{\protect\propositionname}
\theoremstyle{remark}
\newtheorem{rem}[thm]{\protect\remarkname}
\theoremstyle{plain}
\newtheorem{lem}[thm]{\protect\lemmaname}

\usepackage{babel}
\usepackage{mathrsfs}

\usepackage{babel}
\usepackage{mathrsfs}
\usepackage{bbm}

\usepackage{pdfsync}
\numberwithin{equation}{section}

\makeatother

\providecommand{\examplename}{Example}
\providecommand{\lemmaname}{Lemma}
\providecommand{\propositionname}{Proposition}
\providecommand{\remarkname}{Remark}
\providecommand{\theoremname}{Theorem}

\begin{document}
\title{Asymptotic Behavior of the Subordinated Traveling Waves}
\author{\textbf{Yuri Kondratiev}\\
 Department of Mathematics, University of Bielefeld, \\
 D-33615 Bielefeld, Germany,\\
 Dragomanov University, Kiev, Ukraine\\
 Email: kondrat@mathematik.uni-bielefeld.de\\
 Email: kondrat@math.uni-bielefeld.de\and\textbf{ Jos{\'e} Lu{\'i}s
da Silva},\\
 CIMA, University of Madeira, Campus da Penteada,\\
 9020-105 Funchal, Portugal.\\
 Email: joses@staff.uma.pt}
\date{\today}

\maketitle

\begin{abstract}
In this paper we investigate the long-time behavior of the subordination
of the constant speed traveling waves by a general class of kernels.
We use the Feller--Karamata Tauberian theorem in order to study the  long-time 
behavior of the upper and lower wave. As a result we obtain the long-time 
behavior for the propagation of the front of the wave.

\end{abstract}

\tableofcontents{}

\section{Introduction}
\subsection{Object of Study}
Traveling waves form a class of functions which are solutions for different types 
of equations. We have in mind in particular the fractional kinetic corresponding 
to the initial interacting particle system of the Bolker-Pacala model in ecology, 
see \cite{FKKK15} and references therein for more details.  
The present paper is dedicated to study the long-time (or asymptotic) behavior of the propagation of 
the front of the subordinated travelling waves. By a subordination of a solution 
$u(x,t)$ by a density function $G_t(\tau)$, $t,\tau>0$ we mean the function $v(x,t)$ 
defined by  
\[
v(x,t):=\int_0^\infty u(x,\tau)G_t(\tau)\,\mathrm{d}\tau.
\]
The interpretation of the subordination $v(x,t)$ (also called subordination identity 
or subordination principle) is as follows. If the function $u(x,t)$ satisfies an 
evolution equation (say first order time derivative) then under certain conditions, 
$v(x,t)$ satisfies the same type of evolution equation as $u(x,t)$ with the 
first order time derivative replaced by a fractional time derivative. In particular, 
the subordination principle holds for linear PDEs. The fractional derivative 
appearing as a result of subordination is related to the density function $G_t(\tau)$.
In this paper we study three classes (see (C1), (C2), and (C3) below) leading 
to different type of fractional derivatives. These fractional derivatives were 
widely used in physics for modeling slow relaxation and diffusion processes, 
see for example  \cite{MK1,Metzler94,Mainardi2010}. As a simple example 
consider the equation
\begin{equation}
\label{simple-example}
\big(\mathbb{D}_{t}^{\alpha}u_\lambda\big)(t)=-\lambda u_\lambda(t),\quad t>0,\quad u(0)=1,
\end{equation}
where $0<\alpha<1$ and $\mathbb{D}_t^\alpha$ denotes the Caputo-Dzhrbashyan fractional derivative
, see \eqref{eq:Caputo-derivative} for details. It is well known (see for example \cite{KST2006}) that the solution 
of equation \eqref{simple-example} is given in terms of the Mittag-Leffler function $E_\alpha$, namely
\[
u_\lambda(t)=E_\alpha(-\lambda t^\alpha).
\]
It follows from the properties of the Mittag-Leffler function (see \cite{GKMS2014}) that 
$u_\lambda(t)\sim Ct^{-\alpha}$ as $t\to\infty,\;C>0$. Here the simbol $\sim$ means that if $f\sim g$ as 
$t\to\infty$, then  $\lim_{t\to\infty}\frac{f(t)}{g(t)}=1$. In addition, there is a density function $G_t^\alpha(\tau)$ 
such that $u_\lambda(t)$ is a subordination, more precisely
\[
u_\lambda(t)=\int_0^\infty e^{-\lambda\tau}G_t^\alpha(\tau)\,\mathrm{d}\tau,
\]
see Proposition~\ref{exa:distr-alphastab-E} below for more details of $G_t^\alpha(\tau)$. Note that if we replace $\mathbb{D}_{t}^{\alpha}$ by $\frac{\mathrm{d}}{\mathrm{d}t}$ in equation \eqref{simple-example}, then $e^{-\lambda t}$ is the 
 solution of that equation.

\subsection{Description of the Results}
A monotone traveling wave $u(x,t)$ with velocity $v$ in given by a profile function $\psi:\mathbb{R}\longrightarrow[0,1]$ as $u(x,t)=\psi(x-vt)$, $t\ge0$. Without lost of generality we assume that the profile function $\psi$ satisfies 
\[ 
\lim_{t\to -\infty}\psi(t)=1\qquad \mathrm{and}\qquad \lim_{t\to\infty}\psi(t)=0.
\]
For each $\varepsilon>0$ there exist $x_\varepsilon^-,x_\varepsilon^+\in\mathbb{R}$ such that 
\[
u(x,t)<\varepsilon,\;\forall x>x_\varepsilon^+\quad \mathrm{and}\quad u(x,t)>1-\varepsilon,\;\forall x<x_\varepsilon^-.
\] 
This allow us to obtain a lower wave $u_\varepsilon^-(x,t)$ and upper wave $u_\varepsilon^+(x,t)$ such 
that the following chain of inequalities hold
\[
u_\varepsilon^-(x,t)\le u(x,t) \le u_\varepsilon^+(x,t).
\]
 Both the lower and upper wave have an explicitly expression, see Section~\ref{sec:LTB-STW} 
 below for details. Hence, we obtain the chain of inequalities  for the subordination
\[
u_\varepsilon^{E,-}(x,t)\le u^E(x,t) \le u_\varepsilon^{E,+}(x,t).
\] 
The subordination is given with respect to the density of the inverse $E$ of a subordinator. 
\subsection{Motivation: Fractional Kinetic}
One particular way to obtain kinetic equations for densities is the
following, see e.g., \cite{FKK10} for details. Let us consider a
Markov stochastic dynamics for a continuous interacting particle system
in $\mathbb{R}^{d}$. The state evolution of this system may be described
by means of a hierarchical system of evolution equations for correlation
functions. In a mesoscopic scaling limit (e.g., in Vlasov type scaling)
we arrive in the so-called kinetic hierarchy for correlation functions.
Note that, in general, this hierarchy is not related anymore to a Markov
dynamics. But the key property of the kinetic hierarchy is what is
called the chaos preservation in physical literature. In the mathematical
language, it means the following. If we start our system after the
scaling with a Poisson initial measure $\pi_{\rho}$ with the intensity
measure $\mathrm{d}\sigma(x)=\rho(x)\,\mathrm{d}x$, then in the course
of evolutions the state of the system will be again a Poisson measure
$\pi_{\rho_{t}}$ and there exists a non-linear operator $V$ such that the 
density $\rho_{t}(x)$ satisfy the non-linear equation 
\[
\frac{\partial}{\partial t}\rho_{t}(x)=V(\rho_{t})(x),\quad x\in\mathbb{R}^d
\]
with the initial data $\rho_{0}(x)=\rho(x)$. This equation is called
the kinetic equation for the considered stochastic dynamics of an
infinite particle systems. We would like to stress that the kinetic
equation is only one particular byproduct of the kinetic hierarchy
which may be considered as a new important system of equations describing
the dynamics, see comments by H.~Spohn in \cite{Spohn1980}.

Let us consider a random time change in the initial Markov dynamics.
Then we have a hierarchical system of evolution equations with a general
fractional derivative in time corresponding to the random time change,
see \cite{KKS2018,KocKon2017}. After a scaling we obtain a kinetic
hierarchy which is the same as before but with generalized time derivatives
instead of usual ones. This new hierarchy does not preserve anymore the chaos
property. But due to general subordination
principle the solution to the fractional kinetic hierarchy is nothing
but the subordination of the solution to the initial kinetic hierarchy
by a particular kernel associated with the random time. The latter
is deeply related to the linear character of the evolution in the
kinetic hierarchies. As a consequence, we have the evolution of the
density in the fractional dynamics which is nothing but the subordination
of the evolution of the density corresponding to the initial kinetic
equation. Therefore, the kinetic dynamics of the density in the fractional
time is just the transformation of the solution to the kinetic equation
in the physical time. This statement supports our doubts that the
study of non-linear kinetic equations with fractional derivatives
may be justified by arguments coming from physical background. Of
course, we can consider non-linear evolution equations with fractional
derivatives as mathematical objects. But the physical sense of their
solutions remain an open question.

For several particular models of Markov dynamics we already derived
and studied the related kinetic equations, see \cite{FKK10,FKK11,FKKL11}.
In particular, for certain class of such equation we obtained the
existence of solutions in the form of traveling waves, see \cite{Finkelshtein2019a,Finkelshtein2019}.
There appears a natural question about the properties of subordinated
solutions in the case of traveling waves. The physical sense of the
fractional time may be related to a friction included in the initial
system. From the point of view of such interpretation we shall expect
that the motion of the subordinated wave shall be slower comparing
with the initial one. Actually, we will show that this hypothesis
may be justified for several classes of random time changes.

\section{Preliminaries}
\label{sec:preliminaries}In this section we introduce the general
framework we work with. More precisely, we will use the concept of
general fractional derivative (GFD) associated to a kernel $k\in L_{\mathrm{loc}}^{1}(\mathbb{R}_{+})$,
see \cite{Kochubei11}. We consider three classes of admissible kernels
$k$ characterized in terms of their Laplace transforms $\mathcal{K}(\lambda)$
as $\lambda\to0$, see (C1), (C2) and (C3) \vpageref{eq:C1} below.

Let $S=\{S(t),\;t\ge0\}$ be a subordinator without drift, that is,
a process with stationary and independent non-negative increments
starting from $0$, see \cite{Bertoin96} for more details. The Laplace
transform of $S(t)$, $t\ge0$ is expressed as 
\[
\mathbb{E}(e^{-\lambda S(t)})=e^{-t\Phi(\lambda)},\quad\lambda\ge0,
\]
where $\Phi:[0,\infty)\longrightarrow[0,\infty)$ is called the Laplace
exponent which admits the representation 
\begin{equation}
\Phi(\lambda)=\int_{(0,\infty)}(1-e^{-\lambda\tau})\,d\sigma(\tau).\label{eq:Levy-Khintchine}
\end{equation}
The measure $\sigma$ is called L{\'e}vy measure, has support in
$[0,\infty)$ and fulfills 
\begin{equation}
\int_{(0,\infty)}(1\wedge\tau)\,\mathrm{d}\sigma(\tau)<\infty.\label{eq:Levy-condition}
\end{equation}
In what follows we assume that the L{\'e}vy measure $\sigma$ satisfy
\begin{equation}
\sigma(0,\infty)=\infty.\label{eq:Levy-massumption}
\end{equation}
Given the L{\'e}vy measure $\sigma$, we define the function $k$
by 
\begin{equation}
k:(0,\infty)\longrightarrow(0,\infty),\;t\mapsto k(t):=\sigma\big((t,\infty)\big)\label{eq:k}
\end{equation}
and denote its Laplace transform by $\mathcal{K}$, that is, for any
$\lambda\ge0$ one has 
\begin{equation}
\mathcal{K}(\lambda):=\int_{0}^{\infty}e^{-\lambda t}k(t)\,\mathrm{d}t.\label{eq:LT-k}
\end{equation}
The function $\mathcal{K}$ is expressed in terms of the Laplace exponent
$\Phi$ as 
\begin{equation}
\Phi(\lambda)=\lambda\mathcal{K}(\lambda),\quad\forall\lambda > 0.\label{eq:Laplace-exponent}
\end{equation}

\begin{example}
\label{exa:alpha-stable1}
\begin{enumerate}
\item The classical example of a subordinator $S$ is the so-called $\alpha$-stable
process $\alpha\in(0,1)$ with Laplace exponent $\Phi(\lambda)=\lambda^{\alpha}$
and L{\'e}vy measure $d\sigma(\tau)=\frac{\alpha}{\Gamma(1-\alpha)}\tau^{-1-\alpha}\,d\tau.$ 
\item The Gamma process $Y^{(a,b)}$ with parameters $a,b>0$ is another
example of a subordinator with Laplace exponent $\Phi_{(a,b)}(\lambda)=a\log\left(1+\frac{\lambda}{b}\right)$
and L{\'e}vy measure $d\sigma(\tau)=a\tau^{-1}e^{-b\tau}\,d\tau.$ 
\end{enumerate}
\end{example}
Let $E$ be the inverse process of the subordinator $S$, that is,
\begin{equation}
E(t):=\inf\{s\ge0:\;S(s)\ge t\}=\sup\{s\ge0:\;S(s)\le t\}.\label{eq:inverse-sub}
\end{equation}
For any $t\ge0$ we denote by $G_{t}^{k}(\tau)\equiv G_{t}(\tau)$,
$\tau\ge0$ the marginal density of $E(t)$ or, equivalently 
\[
G_{t}(\tau)\,\mathrm{d}\tau=\partial_{\tau}P(E(t)\le\tau)=\partial_{\tau}P(S(\tau)\ge t)=-\partial_{\tau}P(S(\tau)<t).
\]

As the density $G_{t}(\tau)$ plays an important role in the analysis
below here we collect some important properties. 
\begin{prop}[cf.~Prop.~1(a) in \cite{Bingham1971}]
\label{exa:distr-alphastab-E}If $S$ is the $\alpha$-stable process,
$\alpha\in(0,1)$, then the inverse process $E(t)$ has a Mittag-Leffler
distribution, namely 
\begin{equation}
\mathbb{E}(e^{-\lambda E(t)})=\int_{0}^{\infty}e^{-t\tau}G_{t}(\tau)\,\mathrm{d}\tau=\sum_{n=0}^{\infty}\frac{(-\lambda t^{\alpha})^{n}}{\Gamma(n\alpha+1)}=E_{\alpha}(-\lambda t^{\alpha}).\label{eq:Laplace-density-alpha}
\end{equation}
Here $E_\alpha$ is the Mittag-Leffler function with index $\alpha$, see \cite{GKMS2014}.
\end{prop}
\begin{rem}
\label{rem:ML-function} 

\begin{enumerate}
\item It follows from the asymptotic behavior of the Mittag-Leffler function
$E_{\alpha}$ that $\mathbb{E}(e^{-\lambda E(t)})\sim Ct^{-\alpha}$
as $t\to\infty$. 
\item It follows from the properties of the Mittag-Leffler function $E_{\alpha}$,
 that the density $G_{t}(\tau)$ is given in terms of the
Wright function $W_{\mu,\nu}$, namely $G_{t}(\tau)=t^{-\alpha}W_{-\alpha,1-\alpha}(\tau t^{-\alpha})$,
see \cite{Gorenflo1999} for more details. 
\end{enumerate}
\end{rem}
For a general subordinator, the following lemma determines the $t$-Laplace
transform of $G_{t}(\tau)$, with $k$ and $\mathcal{K}$ given in
\eqref{eq:k} and \eqref{eq:LT-k}, respectively. For the proof see
\cite{Kochubei11}. 
\begin{lem}
\label{lem:t-LT-G}The $t$-Laplace transform of the density $G_{t}(\tau)$
is given by 
\begin{equation}
\int_{0}^{\infty}e^{-\lambda t}G_{t}(\tau)\,\mathrm{d}t=\mathcal{K}(\lambda)e^{-\tau\lambda\mathcal{K}(\lambda)}.\label{eq:LT-G-t}
\end{equation}
The double ($\tau,t$)-Laplace transform of $G_{t}(\tau)$ is 
\begin{equation}
\int_{0}^{\infty}\int_{0}^{\infty}e^{-p\tau}e^{-\lambda t}G_{t}(\tau)\,\mathrm{d}t\,\mathrm{d}\tau=\frac{\mathcal{K}(\lambda)}{\lambda\mathcal{K}(\lambda)+p}.\label{eq:double-Laplace}
\end{equation}
\end{lem}
For any $\alpha\in(0,1)$ the Caputo-Dzhrbashyan fractional derivative
of order $\alpha$ of a function $u$ is defined by (see e.g., \cite{KST2006}
and references therein) 
\begin{equation}
\big(\mathbb{D}_{t}^{\alpha}u\big)(t)=\frac{d}{dt}\int_{0}^{t}k(t-\tau)u(\tau)\,\mathrm{d}\tau-k(t)u(0),\quad t>0,\label{eq:Caputo-derivative}
\end{equation}
where 
\[
k(t)=\frac{t^{-\alpha}}{\Gamma(1-\alpha)},\;t>0.
\]
More generally, we consider differential-convolution operators 
\begin{equation}
\big(\mathbb{D}_{t}^{(k)}u\big)(t)=\frac{d}{dt}\int_{0}^{t}k(t-\tau)u(\tau)\,\mathrm{d}\tau-k(t)u(0),\;t>0,\label{eq:general-derivative}
\end{equation}
where $k\in L_{\mathrm{loc}}^{1}(\mathbb{R}_{+})$ is a nonnegative
kernel. The distributed order derivative $\mathbb{D}_{t}^{(\mu)}$
is an example of such operator, corresponding to 
\begin{equation}
k(t)=\int_{0}^{1}\frac{t^{-\tau}}{\Gamma(1-\tau)}\mu(\tau)\,\mathrm{d}\tau,\quad t>0,\label{eq:distributed-kernel}
\end{equation}
where $\mu(\tau)$, $0\le\tau\le1$ is a positive weight function
on $[0,1]$, see \cite{Atanackovic2009,Daftardar-Gejji2008,Gorenflo2005,Hanyga2007,Kochubei2008,Meerschaert2006}.

From now on $L$ always denotes a slowly varying function (SVF) at
infinity (see for instance \cite{Bingham1987} and \cite{Schilling12}),
that is, 
\[
\lim_{x\to\infty}\frac{L(\lambda x)}{L(x)}=1,\qquad\mathrm{for\;any\;}\lambda>0,
\]
while $C$, $C_{\pm}$ are constants whose values are unimportant,
and which may change from line to line.

In the following we consider three classes of admissible kernels $k\in L_{\mathrm{loc}}^{1}(\mathbb{R}_{+})$,
characterized in terms of their Laplace transforms $\mathcal{K}(\lambda)$
as $\lambda\to0$ (i.e., as local conditions): 
\begin{equation}
\mathcal{K}(\lambda)\sim\lambda^{\alpha-1},\quad0<\alpha<1.\tag*{(C1)}\label{eq:C1}
\end{equation}
\begin{equation}
\mathcal{K}(\lambda)\sim\lambda^{-1}L\left(\frac{1}{\lambda}\right),\quad L(y):=\mu(0)\log(y)^{-1}.\tag*{(C2)}\label{eq:C2}
\end{equation}
\begin{equation}
\mathcal{K}(\lambda)\sim\lambda^{-1}L\left(\frac{1}{\lambda}\right),\quad L(y):=C\log(y)^{-1-s},\;s>0,\;C>0.\tag*{(C3)}\label{eq:C3}
\end{equation}
We would like to emphasize that these three classes of kernels lead
to different differential-convolution operators. In particular, the
Caputo-Djrbashian fractional derivative (C1) and distributed order
derivatives (C2), (C3). In the next section we study the long-time
behavior of the subordination of the constant speed traveling wave 
corresponding to these differential-convolution
operators. Working in such generality a price must be paid, namely
the replacement of the fundamental solution by its Cesaro mean. This
is the key technical observation that underlies the analysis of several
different model situations.

\section{Long-time Behavior of the Subordination of Traveling Waves}
\label{sec:LTB-STW}
A monotone traveling wave $u(x,t)$, $t\ge0$, $x\in\mathbb{R}$ with
velocity $v>0$ is defined by a profile function $\psi:\mathbb{R}\longrightarrow[0,1]$
which is a continuous monotonically decreasing function such that 
\[
\lim_{x\to-\infty}\psi(x)=1, \qquad \lim_{x\to\infty}\psi(x)=0,
\]
and $u(x,t)=\psi(x-vt)$, $t\ge0$ for almost all $x\in\mathbb{R}$.
For each $\varepsilon>0$ and $t>0$ introduce $x_{\varepsilon}^{-},x_{\varepsilon}^{+}\in\mathbb{R}$
as 
\[
\forall x>x_{\varepsilon}^{+},\;\;u(x,t)<\varepsilon\quad\mathrm{and}\quad\forall x<x_{\varepsilon}^{-},\;\;u(x,t)>1-\varepsilon
\]
and consider the two-side estimate of $u(x,t)$ for any $x\in\mathbb{R}$
and $t\ge0$ 
\begin{equation}
u_{\varepsilon}^{-}(x,t)\le u(x,t)\le u_{\varepsilon}^{+}(x,t),\label{eq:TV-two-side-estimate}
\end{equation}
where 
\begin{eqnarray*}
u_{\varepsilon}^{+}(x,t)&:=&\mathbbm{1}_{(-\infty,x_{\varepsilon}^{+}]}(x-vt)+\varepsilon\mathbbm{1}_{[x_{\varepsilon}^{+},\infty)}(x-vt),\\
u_{\varepsilon}^{-}(x,t)&:=&(1-\varepsilon)\mathbbm{1}_{(-\infty,x_{\varepsilon}^{-}]}(x-vt).
\end{eqnarray*}

\begin{figure}
\begin{centering}
\includegraphics{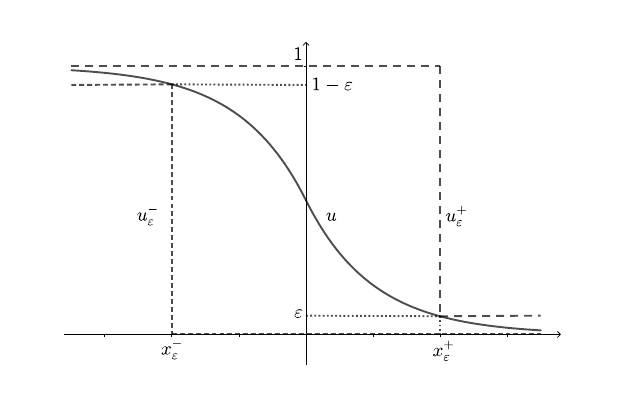}
\par\end{centering}
\caption{\label{fig:Traveling-wave}Traveling wave $u$, upper wave $u_{\varepsilon}^{+}$
and lower wave $u_{\varepsilon}^{-}$.}
\end{figure}

The functions $u_{\varepsilon}^{-}(x,t)$ and $u_{\varepsilon}^{+}(x,t)$
we will call lower and upper waves, respectively, see Figure~\ref{fig:Traveling-wave}.

\subsection{Subordination of the Traveling Wave}
Consider the solutions of the evolution equations 
\begin{eqnarray}
\frac{\partial}{\partial t}u_1(x,t)&=&(Au_1)(x,t) \label{evolution1}\\
\big(\mathbb{D}_{t}^{(k)}u_k\big)(t)&=&(Au_k)(x,t),\label{evolution2}
\end{eqnarray}
where $A$ is an operator acting in the spatial variable $x$ and the same initial conditions
\[
u_1(x,0)=\xi(x),\qquad u_k(x,0)=\xi(x).
\] 
Then in certain conditions (e.g.~$A$ closed linear operator) the solutions $u_1$ and $u_k$ 
satisfy the subordination identity (also known as subordination principle), that is, there 
exists a nonnegative function $G_t(\tau)$, $t,\tau>0$ such that $\int_0^\infty G_t(\tau)\,\mathrm{d}\tau=1$ 
and 
\[
u_k(x,t)=\int_0^\infty u_1(x,\tau)G_t(\tau)\,\mathrm{d}\tau.
\] 
aThe proper notion of solution for equations \eqref{evolution1} and \eqref{evolution2} were explained 
in \cite{Kochubei11} when $A$ is the Laplace operator on $\mathbb{R}^n$, in \cite{Bazhlekova00,Baz01,Bazhlekova2015} in the framework of semigroups generators (for special classes of $k$) or in \cite{Pruss12} for abstract Volterra equations. 

We are interested in the subordination of the traveling wave $u(x,t)$ by the density 
$G_t(\tau)$ associated to the inverse process $E$. Hence, we obtain a new function $u^E(x,t)$ defined by 
\begin{equation}
u^{E}(x,t):=\int_{0}^{\infty}u(x,\tau)G_{t}(\tau)\,\mathrm{d}\tau.\label{eq:sub-TV}
\end{equation}
The subordination of the lower and upper waves, denoted by $u^{E,-}_\varepsilon(x,t)$ 
and  $u^{E,+}_\varepsilon(x,t)$, respectively,  are defined similarly. Having in mind the chain of inequalities 
\eqref{eq:TV-two-side-estimate} we obtain the chain for the subordinated functions
\begin{equation}
u^{E,-}_\varepsilon(x,t) \le u^{E}(x,t) \le u^{E,+}_\varepsilon(x,t).
\end{equation}

\begin{rem}
\label{rem:two-estimates}
The long time behavior of the 
function $u^E(x,t)$ as  $t\to \infty$ may be determined, 
under certain conditions, by studying the behavior of its Laplace transform
 $\tilde{u}^E(x,\lambda)$ as $\lambda\to 0$, and vice versa. 
An important situation where such a correspondence holds is described by the 
 Feller--Karamata Tauberian (FKT) theorem. 
 \end{rem}
 
 We state below a version of the FKT 
 theorem  which suffices for our purposes, see the monographs \cite[Sec.~1.7]{Bingham1987}
and \cite[XIII, Sec.~1.5]{Feller71} for a more general version and proofs.

\begin{thm}[Feller--Karamata Tauberian]
\label{thm:FKT-LST}Let $U:[0,\infty)\longrightarrow\mathbb{R}$
be a monotone non-decreasing right-continuous function such that 
\[
w(\lambda):=\int_{0}^{\infty}e^{-\lambda t}\,\mathrm{d}U(t)<\infty,\quad\forall\lambda>0.
\]
If $L$ is a slowly varying function and  $C,\rho\ge0$, then the following
are equivalent
\begin{equation}
U(t)\sim\frac{C}{\Gamma(\rho+1)}t^{\rho}L(t)\quad\mathrm{as}\;t\to\infty,\label{eq:asymp-U}
\end{equation}
\begin{equation}
w(\lambda)\sim C\lambda^{-\rho}L\left(\frac{1}{\lambda}\right)\quad\mathrm{as}\;\lambda\to0^{+}.\label{eq:asym-w}
\end{equation}
When $C=0$, (\ref{eq:asymp-U}) is to be interpreted as $U(t)=o(t^{\rho}L(t))$;
similarly for (\ref{eq:asym-w}). 
\end{thm}

\begin{rem}
\label{rem:Subord-Twave}
\begin{enumerate}
\item In general, the function $u^{E}(x,t)$ is not monotone in $t$,
that will be needed to apply the  Theorem~\ref{thm:FKT-LST}. In addition, 
the Laplace transform of $u^{E}(x,t)$ can be explicitly computed only when $G_t(\tau)$
corresponds to the density of the inverse stable subordinator. 
\item We define the $t$-increasing function 
\[
\int_{0}^{t}u^{E}(x,s)\,\mathrm{d}s
\]
and then will obtain the long-time behavior for the Cesaro mean of
$u^{E}(x,t)$, that is, 
\begin{equation}
\label{Cesaro-mean}
M_{t}(u^{E}(x,\cdot)):=\frac{1}{t}\int_{0}^{t}u^{E}(x,s)\,\mathrm{d}s.
\end{equation}
\item For any fixed time $t$, the subordinated traveling wave $u^{E}(x,t)$
is decreasing and continuous in $x$. Hence, given $\beta\in(0,1)$
there is a unique $x_{\beta}^{E}(t)\in(0,1)$ which solves the equation
\[
u^{E}(x_{\beta}^{E}(t),t)=\beta.
\]
We call $x_{\beta}^{E}(t)$ the propagation of the front of $u^{E}(x,t)$
of the level $\beta$. For a general definition of the \emph{propagation
of the front} of a function $u(x,t)$, which is the solution of a certain
differential equation, see \cite{FKT2019} and references therein. 
\end{enumerate}
\end{rem}

The considerations in Remark~\ref{rem:Subord-Twave} lead us 
to consider the chain of inequalities for the Cesaro means, namely
\[
M_{t}(u^{E,-}_\varepsilon(x,\cdot))\le M_{t}(u^{E}(x,\cdot)) \le M_{t}(u^{E,+}_\varepsilon(x,\cdot)).
\]
Unfortunately the FKT theorem does not apply to inequalities. Hence,
we study the long-time behavior of the Cesaro mean of the subordination of 
both the upper and lower waves separately. It turns out that both of these long-time 
behavior are of the same type, compare for example \eqref{eq:CM-long-t-lb} 
and \eqref{eq:CM-long-t-ub} for the class (C1). Although we are not allowed to conclude any type of 
long-time behavior for the Cesaro mean of our subordination travelling wave $u^E(x,t)$, it gives 
good indications and we may derive a two-side estimation for the propagation of the front 
which are again of the type. 
The results for the class (C1) are stated in Theorem~\ref{thm:FPropagation-C1}
below, see also Theorem~\ref{thm:FPropagation-C2} (resp.~Theorem~\ref{thm:FPropagation-C3}) 
for the class (C2) (resp.~class (C3)).

\subsection{Long-Time Behavior: Class (C1)}

\subsubsection{The Subordination of the Lower Wave}

We start with the lower wave, namely the subordination 
\begin{equation} 
\label{eq:sub-lower-bound}
u_{\varepsilon}^{E,-}(x,t):=(1-\varepsilon)\int_{0}^{\infty}\mathbbm{1}_{(-\infty,x_{\varepsilon}^{-}]}(x-v\tau)G_{t}(\tau)\,\mathrm{d}\tau.
\end{equation}
If we denote by $\theta_{\varepsilon}^{-}:=\frac{x-x_{\varepsilon}^{-}}{v}$
with $x>x_{\varepsilon}^{-}$, then $u_{\varepsilon}^{E,-}(x,t)$
is given by 
\[
u_{\varepsilon}^{E,-}(x,t)=(1-\varepsilon)\int_{\theta_{\varepsilon}^{-}}^{\infty}G_{t}(\tau)\,\mathrm{d}\tau.
\]
Computing the Laplace transform of the monotone function $v^{-}(x,t):=\int_{0}^{t}u_{\varepsilon}^{E,-}(x,s)\,\mathrm{d}s$
\begin{eqnarray*}
(\mathscr{L}v^-(x,\cdot))(\lambda)&=&\int_0^\infty e^{-\lambda t}\,\mathrm{d}_t v^-(x,t)=\int_0^\infty e^{-\lambda t}u^{E,-}_\varepsilon(x,t)\,\mathrm{d}t  \\ 
&=&\int_0^\infty e^{-\lambda t}\int_{\theta_\varepsilon^-}^\infty G_t(\tau)\,\mathrm{d}\tau\,\mathrm{d}t. 
\end{eqnarray*}
Using Fubini theorem and equality \eqref{eq:LT-G-t} yields
\begin{equation}
(\mathscr{L}v^-(x,\cdot))(\lambda)=(1-\varepsilon)\mathcal{K}(\lambda)\int_{\theta_{\varepsilon}^{-}}^{\infty}e^{-\tau\lambda\mathcal{K}(\lambda)}\,\mathrm{d}\tau=(1-\varepsilon)\lambda^{-1}e^{-\theta_{\varepsilon}^{-}\lambda\mathcal{K}(\lambda)}.\label{eq:LT-lower-bound}
\end{equation}

For the class (C1) we have $\mathcal{K}(\lambda)\sim \lambda^{\alpha-1}$, $\lambda\to 0$,
$0<\alpha<1$, hence 
\[
(\mathscr{L}v^-(x,\cdot))(\lambda)\sim (1-\varepsilon)\lambda^{-1}e^{-\theta_{\varepsilon}^{-}\lambda^{\alpha}}=\lambda^{-\rho}L\left(\frac{1}{\lambda}\right),\quad \lambda\to 0,
\]
where $\rho=1$ and $L(y)=(1-\varepsilon)\exp(-\theta_{\varepsilon}^{-}y^{-\alpha})$
is a SVF. Then we conclude by FKT (see Theorem \ref{thm:FKT-LST}) that $v^-({x,t})\sim tL(t)$, $t\to\infty$ which implies the long-time behavior for the Cesaro mean of the subordination of the lower wave 
$u^{E,-}_\varepsilon(x,t)$:
\begin{equation}
M_{t}(u_{\varepsilon}^{E,-}(x,\cdot))=\frac{1}{t}v^-(x,t)\sim L(t)=(1-\varepsilon)e^{-\theta_\varepsilon^- t^{-\alpha}},\;t\to\infty.\label{eq:CM-long-t-lb}
\end{equation}
Define the right hand side of the above by $W_{\varepsilon}^{-}(x,t)$,
that is, 
\[
W_{\varepsilon}^{-}(x,t):=(1-\varepsilon)e^{-\theta_\varepsilon^- t^{-\alpha}}.
\]
It is clear that for any fixed $x$ we have 
\[
W_{\varepsilon}^{-}(x,t)\to1-\varepsilon,\quad t\to\infty
\]
and fixing $t$ (recall $\theta_\varepsilon^-=\frac{x-x_{\varepsilon}^{-}}{v}$) yields 
\[
W_{\varepsilon}^{-}(x,t)\to0,\quad x\to\infty.
\]
To find the propagation of the front $x_{\varepsilon,\beta}^{-}(t)$
of the Cesaro mean $M_{t}(u_{\varepsilon}^{E,-}(x,\cdot))$ of the
level $\beta\in(0,1-\varepsilon)$ we solve the equation 
\[
W_{\varepsilon}^{-}(x_{\varepsilon,\beta}^{-}(t),t)=\beta
\]
to obtain 
\[
x_{\varepsilon,\beta}^{-}(t)\sim vt^{\alpha}\log\left(\frac{1-\varepsilon}{\beta}\right)+x_{\varepsilon}^{-}=:C_{-}t^{\alpha}+x_{\varepsilon}^{-},\quad t\to\infty.
\]
So, the propagation of the front of the Cesaro mean $M_{t}(u_{\varepsilon}^{E,-}(x,\cdot))$
 is $x_{\varepsilon,\beta}^{-}(t)\sim C_{-}t^{\alpha}$
as $t\to\infty$.

\subsubsection{The Subordination of the Upper Wave}

We are now interested in the upper wave, namely the subordination
\begin{align}
u_{\varepsilon}^{E,+}(x,t) & :=\int_{0}^{\infty}u^+_\varepsilon(x,\tau)G_{t}(\tau)\,\mathrm{d}\tau\label{eq:sub-upper-bound}\\
 & =\int_{0}^{\infty}\mathbbm{1}_{(-\infty,x_{\varepsilon}^{+}]}(x-v\tau)G_{t}(\tau)\,\mathrm{d}\tau+\varepsilon\int_{0}^{\infty}\mathbbm{1}_{[x_{\varepsilon}^{+},\infty)}(x-v\tau)G_{t}(\tau)\,\mathrm{d}\tau\\
 & =:u_{\varepsilon}^{E,+,1}(x,t)+u_{\varepsilon}^{E,+,2}(x,t).
\end{align}
As before we study the Cesaro mean of each of the above functions,
namely 
\begin{align*}
M_{t}(u_{\varepsilon}^{E,+,1}(x,\cdot)) & :=\frac{1}{t}\int_{0}^{t}u_{\varepsilon}^{E,+,1}(x,s)\,\mathrm{d}s,\\
M_{t}(u_{\varepsilon}^{E,+,2}(x,\cdot)) & :=\frac{1}{t}\int_{0}^{t}u_{\varepsilon}^{E,+,2}(x,s)\,\mathrm{d}s.
\end{align*}
If we denote by $\theta_{\varepsilon}^{+}:=\frac{x-x_{\varepsilon}^{+}}{v}$
with $x>x_{\varepsilon}^{+}$, then $u_{\varepsilon}^{E,+,1}(x,t)$
is equal to 
\[
u_{\varepsilon}^{E,+,1}(x,t):=\int_{\theta_{\varepsilon}^{+}}^{\infty}G_{t}(\tau)\,\mathrm{d}\tau.
\]
Computing the Laplace transform of the monotone function $v^{+,1}(x,t):=\int_{0}^{t}u_{\varepsilon}^{E,+,1}(x,s)\,\mathrm{d}s$
and using \eqref{eq:LT-G-t} yields 
\begin{equation}
(\mathscr{L}v^{+,1}(x,\cdot))(\lambda)=\mathcal{K}(\lambda)\int_{\theta_{\varepsilon}^{+}}^{\infty}e^{-\tau\lambda\mathcal{K}(\lambda)}\,\mathrm{d}\tau=\lambda^{-1}e^{-\theta_{\varepsilon}^{+}\lambda\mathcal{K}(\lambda)}.\label{eq:LT-upper-bound-1}
\end{equation}
A similar procedure for $v^{+,2}(x,t):=\int_{0}^{t}u_{\varepsilon}^{E,+,2}(x,s)\,\mathrm{d}s$
produces the following Laplace transform 
\begin{equation}
(\mathscr{L}v^{+,2}(x,\cdot))(\lambda)=\varepsilon\mathcal{K}(\lambda)\int_{0}^{\theta_{\varepsilon}^{+}}e^{-\tau\lambda\mathcal{K}(\lambda)}\,\mathrm{d}\tau=\varepsilon\lambda^{-1}(1-e^{-\theta_{\varepsilon}^{+}\lambda\mathcal{K}(\lambda)}).\label{eq:LT-upper-bound-2}
\end{equation}

For the class (C1), $\mathcal{K}(\lambda)\sim\lambda^{\alpha-1}$, $\lambda\to 0$, $0<\alpha<1$,
it follows from \eqref{eq:LT-upper-bound-1} that 
\[
(\mathscr{L}u_{\varepsilon}^{E,+,1}(x,\cdot))(\lambda)\sim \lambda^{-1}e^{-\theta_{\varepsilon}^{+}\lambda^{\alpha}}=\lambda^{-\rho}L\left(\frac{1}{\lambda}\right),
\]
where $\rho=1$ and $L(y)=\exp(-\theta_{\varepsilon}^{+}y^{-\alpha})$
is a SVF. Then we conclude 
\begin{equation}
M_{t}(u_{\varepsilon}^{E,+,1}(x,\cdot))\sim L(t)=e^{-\theta_{\varepsilon}^{+}t^{-\alpha}},\;t\to\infty.\label{eq:CM-long-t-ub-1}
\end{equation}
For the second function $u_{\varepsilon}^{E,+,2}(x,t)$ we obtain
\[
(\mathscr{L}u_{\varepsilon}^{E,+,2}(x,\cdot))(\lambda)=\varepsilon\lambda^{-1}(1-e^{-\theta_{\varepsilon}^{+}\lambda^{\alpha}})=\lambda^{-\rho}L\left(\frac{1}{\lambda}\right),
\]
where $\rho=1$ and $L(y)=\varepsilon(1-\exp(-\theta_{\varepsilon}^{+}y^{-\alpha})$
is a SVF. Thus 
\begin{equation}
M_{t}(u_{\varepsilon}^{E,+,2}(x,\cdot))\sim L(t)=\varepsilon\left(1-e^{-\theta_{\varepsilon}^{+}t^{-\alpha}}\right),\;t\to\infty.\label{eq:CM-long-t-ub-2}
\end{equation}
Putting \eqref{eq:CM-long-t-ub-1} and \eqref{eq:CM-long-t-ub-2}
together we obtain 
\begin{equation}
\label{eq:CM-long-t-ub}
M_{t}(u_{\varepsilon}^{E,+}(x,\cdot))\sim(1-\varepsilon)\exp\left(-\frac{x-x_{\varepsilon}^{+}}{v}t^{-\alpha}\right)+\varepsilon.
\end{equation}
Define the right hand side of the above by $W_{\varepsilon}^{+}(x_{\varepsilon}^{+},t)$,
that is 
\[
W_{\varepsilon}^{+}(x_{\varepsilon}^{+},t):=(1-\varepsilon)\exp\left(-\frac{x-x_{\varepsilon}^{+}}{v}t^{-\alpha}\right)+\varepsilon.
\]
For any fixed $x$ we have 
\[
W_{\varepsilon}^{+}(x_{\varepsilon}^{+},t)\to1,\quad t\to\infty
\]
and if $t$ is fixed we obtain 
\[
W_{\varepsilon}^{+}(x_{\varepsilon}^{+},t)\to\varepsilon,\quad x\to\infty.
\]
To find the propagation of the front $x_{\varepsilon,\beta}^{+}(t)$
of the Cesaro mean $M_{t}(u_{\varepsilon}^{E,+}(x,\cdot))$ of the
level $\beta\in(\varepsilon,1)$ we solve the equation 
\[
W_{\varepsilon}^{+}(x_{\varepsilon,\beta}^{+}(t),t)=\beta
\]
for $x_{\varepsilon,\beta}^{+}(t)$ and obtain 
\[
x_{\varepsilon,\beta}^{+}(t)=vt^{\alpha}\log\left(\frac{1-\varepsilon}{\beta-1}\right)+x_{\varepsilon}^{+}=:C_{+}t^{\alpha}+x_{\varepsilon}^{-},\quad t\to\infty.
\]
So, the propagation of the front of the Cesaro mean $M_{t}(u_{\varepsilon}^{E,-}(x,t))$
of the subordinated of the upper wave is $x_{\varepsilon,\beta}^{+}(t)\sim C_{+}t^{\alpha}$
as $t\to\infty$. 
\begin{rem}
For any $x\in\mathbb{R}$ and $t\ge0$ we have the following chain
of inequalities, cf.~\eqref{eq:TV-two-side-estimate} 
\[
u_{\varepsilon}^{E,-}(x,t)\le u^{E}(x,t)\le u_{\varepsilon}^{E,+}(x,t).
\]
As $G_{t}(\tau)$ is a density, the same type of chain for the Cesaro
mean is also valid, that is 
\[
M_{t}\big(u_{\varepsilon}^{E,-}(x,\cdot)\big)\le M_{t}\big(u^{E}(x,\cdot)\big)\le M_{t}\big(u_{\varepsilon}^{E,+}(x,\cdot)\big).
\]
If $x_{\beta}(t)\in(\varepsilon,1-\varepsilon)$ denotes the propagation
of the front of the Cesaro mean $M_{t}\big(\psi^{E}(x,t)\big)$ of
the level $\beta$, then the following relation between the propagation
of the fronts of the level $\beta$ hold 
\[
C_{-}t^{\alpha}\sim x_{\varepsilon}^{-}(t)\le x_{\beta}(t)\le x_{\varepsilon,\beta}^{+}(t)\sim C_{+}t^{\alpha},\quad t\to \infty.
\]
\end{rem}
We have shown the following theorem. 
\begin{thm}
\label{thm:FPropagation-C1}Let $u(x,t)=\psi(x-vt)$ be a traveling
wave with constant speed $v$ with two-side estimate, for any $x\in\mathbb{R}$,
$t\ge0$ and $\varepsilon>0$ 
\[
u_{\varepsilon}^{-}(x,t)\le u(x,t)\le u_{\varepsilon}^{+}(x,t).
\]
The subordination $u^{E}(x,t)$ of $u(x,t)$ with the density
$G_{t}(\tau)$ (corresponding to the class (C1)) has Cesaro mean $M_{t}\big(u^{E}(x,\cdot)\big)$
with a propagation of the front $x_{\beta}(t)$ of the level $\beta$
that satisfies the two-side estimate, with $C_{-},C_{+}>0$
\[
C_{-}t^{\alpha}\le x_{\beta}(t)\le C_{+}t^{\alpha}\quad \mathrm{as}\quad t\to\infty.
\]
Then we have $x_{\beta}(t)\sim Ct^{\alpha}$ as $t\to\infty$ with $C>0$.
\end{thm}
\begin{rem}
\label{rem:class-C1}In Example 5 of \cite{DKT2018} it is shown, using
a direct method, for the particular example of the inverse stable
subordinator from the class (C1) the long-time behavior of the propagation
of the front $x_{\beta}(t)$ is given by 
\[
x_{\beta}(t)=Ct^{\alpha}+o(t^{\alpha}),\quad t\to\infty.
\]
This shows that when the long-time behavior exists for the subordinated
wave, then the Cesaro mean gives the right result. 
\end{rem}

\subsection{Long-Time Behavior: Class (C2)}

\subsubsection{The Subordination of the Lower Wave}

Here we have $\mathcal{K}(\lambda)\sim\lambda^{-1}L(\lambda^{-1})$
as $\lambda\to0$, where $L(y)=\mu(0)\log(y)^{-1}$, $\mu(0)\neq0$.
It follows from \eqref{eq:LT-lower-bound} that 
\[
\left(\mathscr{L}\int_{0}^{\cdot}u_{\varepsilon}^{E,-}(x,s)\,\mathrm{d}s\right)(\lambda)=(1-\varepsilon)\lambda^{-1}e^{-\theta_{\varepsilon}^{-}\mu(0)\log(\lambda^{-1})^{-1}}=\lambda^{-\rho}L\left(\frac{1}{\lambda}\right),
\]
where $\rho=1$ and $L(y)=(1-\varepsilon)\exp(-\theta_{\varepsilon}^{-}\mu(0)\log(y)^{-1})$
is a SVF. From this follows 
\[
M_{t}(u_{\varepsilon}^{E,-}(x,\cdot))\sim(1-\varepsilon)\exp\left(-\theta_{\varepsilon}^{-}\mu(0)\log(t)^{-1}\right),\;t\to\infty.
\]
For this class of kernels $k$, the propagation of the front $x_{\varepsilon,\beta}^{-}(t)$
of the Cesaro mean $M_{t}(u_{\varepsilon}^{E,-}(x,t))$ of the
level $\beta\in(0,1-\varepsilon)$ solves 
\[
(1-\varepsilon)\exp\left(-\frac{x_{\varepsilon,\beta}^{-}(t)-x_{\varepsilon}^{-}}{v}\mu(0)\log(t)^{-1}\right)=\beta.
\]
We obtain 
\[
x_{\varepsilon,\beta}^{-}(t)=\log\left(\frac{1-\varepsilon}{\beta}\right)\frac{v}{\mu(0)}\log(t)+x_{\varepsilon}^{-}=C_{-}\log(t)+x_{\varepsilon}^{-}.
\]
from which follows the propagation of the front $x_{\varepsilon,\beta}^{-}(t)\sim C_{-}\log(t)$
as $t\to\infty$.

\subsubsection{The Subordination of the Upper Wave}

We have $\mathcal{K}(\lambda)\sim\lambda^{-1}L(\lambda^{-1})$ as
$\lambda\to0$, where $L(y)=\mu(0)\log(y)^{-1}$, $\mu(0)\neq0$.
It follows from \eqref{eq:LT-upper-bound-1} that 
\[
\left(\mathscr{L}\int_{0}^{\cdot}u_{\varepsilon}^{E,+,1}(x,s)\,\mathrm{d}s\right)(\lambda)=\lambda^{-1}e^{-\theta_{\varepsilon}^{+}\lambda\mathcal{K}(\lambda)}=\lambda^{-1}e^{-\theta_{\varepsilon}^{+}\log(\lambda^{-1})^{-1}}=\lambda^{-\rho}L\left(\frac{1}{\lambda}\right),
\]
where $\rho=1$ and $L(y)=\exp(-\theta_{\varepsilon}^{+}\log(y)^{-1})$
is a SVF. From this follows 
\[
M_{t}(u_{\varepsilon}^{E,+,1}(x,\cdot))\sim\exp\left(-\theta_{\varepsilon}^{+}\log(t)^{-1}\right),\;t\to\infty.
\]
For the function $u_{\varepsilon}^{E,+,2}(x,t)$ we obtain 
\[
\left(\mathscr{L}\int_{0}^{\cdot}u_{\varepsilon}^{E,+,2}(x,s)\,\mathrm{d}s\right)(\lambda)=\varepsilon\lambda^{-1}(1-e^{-\theta_{\varepsilon}^{+}\lambda\mathcal{K}(\lambda)})=\varepsilon\lambda^{-1}(1-e^{-\theta_{\varepsilon}^{+}\log(\lambda^{-1})^{-1}})=\lambda^{-\rho}L\left(\frac{1}{\lambda}\right),
\]
where $\rho=1$ and $L(y)=\varepsilon(1-\exp(-\theta_{\varepsilon}^{+}\log(y)^{-1}))$
is a SVF. Hence, we have 
\[
M_{t}(u_{\varepsilon}^{E,+,2}(x,\cdot))\sim\varepsilon\left(1-\exp\left(-\theta_{\varepsilon}^{+}\log(t)^{-1}\right)\right),\;t\to\infty.
\]
Putting together, we obtain the long-time behavior of the Cesaro mean
of $u_{\varepsilon}^{E,+}(x,t)$ for the class (C2), namely 
\[
M_{t}(u_{\varepsilon}^{E,+}(x,\cdot))\sim(1-\varepsilon)\exp\left(-\frac{x-x_{\varepsilon}^{+}}{v}\log(t)^{-1}\right)+\varepsilon,\;t\to\infty.
\]
For this class of kernels $k$, the propagation of the front $x_{\varepsilon,\beta}^{+}(t)$
of the Cesaro mean $M_{t}(u_{\varepsilon}^{E,+}(x,\cdot))$ of the
level $\beta\in(\varepsilon,1)$ is the solution of 
\[
W_{\varepsilon}^{+}(x_{\varepsilon,\beta}^{+}(t),t)=(1-\varepsilon)\exp\left(-\frac{x_{\varepsilon,\beta}^{+}(t)-x_{\varepsilon}^{+}}{v}\log(t)^{-1}\right)+\varepsilon=\beta.
\]
solving for $x_{\varepsilon}^{+}(t)$ we obtain 
\[
x_{\varepsilon,\beta}^{+}(t)=\log\left(\frac{1-\varepsilon}{\beta-\varepsilon}\right)v\log(t)+x_{\varepsilon}^{+}=C_{+}\log(t)+x_{\varepsilon}^{+}.
\]
from which follows the the propagation of the front $x_{\varepsilon,\beta}^{+}(t)\sim C_{+}\log(t)$
as $t\to\infty$. This agrees with the propagation of the front for
the lower bound.

We summarize the results for the class (C2) in the following theorem. 
\begin{thm}
\label{thm:FPropagation-C2}Let $u(x,t)=\psi(x-vt)$ be a traveling
wave with constant speed $v$ with two-side estimate, for any $x\in\mathbb{R}$,
$t\ge0$ and $\varepsilon>0$ 
\[
u_{\varepsilon}^{-}(x,t)\le u(x,t)\le u_{\varepsilon}^{+}(x,t).
\]
The subordination $u^{E}(x,t)$ of $u(x,t)$ with the density
$G_{t}(\tau)$ (corresponding to the class (C2)) has Cesaro mean $M_{t}\big(u^{E}(x,\cdot)\big)$
has a propagation of the front $x_{\beta}(t)$ of the level $\beta$
that satisfies the two-side estimate, with $C_{-},C_{+}>0$
\[
C_{-}\log(t)\le x_{\beta}(t)\le C_{+}\log(t)\quad \mathrm{as}\quad t\to\infty.
\]
Then we have $x_{\beta}(t)\sim C\log(t)$ as $t\to\infty$ with $C>0$.
\end{thm}

\subsection{Long-Time Behavior: Class (C3)}

\subsubsection{The Subordination of the Lower Wave}

We  have asymptotic $\mathcal{K}(\lambda)\sim C\lambda^{-1}L(\lambda^{-1})^{-1-s}$
as $\lambda\to0$ and $s>0$, $C>0$. The substitution of this $\mathcal{K}(\lambda)$
in \eqref{eq:LT-lower-bound} produces 
\[
\left(\mathscr{L}\int_{0}^{\cdot}u_{\varepsilon}^{E,-}(x,s)\,\mathrm{d}s\right))(\lambda)=(1-\varepsilon)\lambda^{-1}e^{-C\theta_{\varepsilon}^{-}L(\lambda^{-1})^{-1-s}}=\lambda^{-\rho}L\left(\frac{1}{\lambda}\right),
\]
where $\rho=1$ and $L(y)=(1-\varepsilon)\exp(-C\theta_{\varepsilon}^{-}\log(y)^{-1-s})$
is a SVF. We conclude that 
\[
M_{t}(u_{\varepsilon}^{E,-}(x,\cdot))\sim(1-\varepsilon)\exp\left(-C\theta_{\varepsilon}^{-}\log(t)^{-1-s}\right),\;t\to\infty.
\]
To find the propagation of the front $x_{\varepsilon,\beta}^{-}(t)$
of $M_{t}(u_{\varepsilon}^{E,-}(x,\cdot))$ of the level $\beta\in(0,1-\varepsilon)$
we solve the equation 
\[
(1-\varepsilon)\exp\left(-C\frac{x_{\varepsilon,\beta}^{-}(t)-x_{\varepsilon}^{-}}{v}\log(t)^{-1-s}\right)=\beta
\]
and obtain 
\[
x_{\varepsilon,\beta}^{-}(t)=\log\left(\frac{1-\varepsilon}{\beta}\right)\frac{v}{C}\log(t)^{1+s}+x_{\varepsilon}^{-}=:C_{-}\log(t)^{1+s}+x_{\varepsilon}^{-}.
\]
Hence, the propagation of the front is $x_{\varepsilon,\beta}^{-}(t)\sim C_{-}\log(t)^{1+s}$
as $t\to\infty$.

\subsubsection{The Subordination of the Upper Wave}

As $\mathcal{K}(\lambda)\sim C\lambda^{-1}L(\lambda^{-1})^{-1-s}$
as $\lambda\to0$ and $s>0$, $C>0$, the substitution of this $\mathcal{K}(\lambda)$
in \eqref{eq:LT-upper-bound-1} produces 
\[
\left(\mathscr{L}\int_{0}^{\cdot}\psi_{\varepsilon}^{E,+,1}(x,\cdot)\,\mathrm{d}s\right)(\lambda)=\lambda^{-1}e^{-C\theta_{\varepsilon}^{+}L(\lambda^{-1})^{-1-s}}=\lambda^{-\rho}L\left(\frac{1}{\lambda}\right),
\]
where $\rho=1$ and $L(y)=\exp(-C\theta_{\varepsilon}^{+}\log(y)^{-1-s})$
is a SVF. We conclude that 
\[
M_{t}(u_{\varepsilon}^{E,+,1}(x,\cdot))\sim\exp\left(-C\theta_{\varepsilon}^{+}\log(t)^{-1-s}\right),\;t\to\infty.
\]
For the function $u_{\varepsilon}^{E,+,2}(x,t)$ we obtain 
\[
\left(\mathscr{L}\int_{0}^{\cdot}u_{\varepsilon}^{E,+,2}(x,s)\,\mathrm{d}s\right)(\lambda)=\varepsilon\lambda^{-1}(1-e^{-C\theta_{\varepsilon}^{+}L(\lambda^{-1})^{-1-s}}),
\]
where $\rho=1$ and $L(y)=\varepsilon(1-\exp(-C\theta_{\varepsilon}^{+}\log(y)^{-1-s})$
is a SVF. Hence, we conclude that 
\[
M_{t}(u_{\varepsilon}^{E,+,2}(x,\cdot))\sim\varepsilon\left(1-\exp\left(-C\theta_{\varepsilon}^{+}\log(t)^{-1-s}\right)\right),\;t\to\infty.
\]
Therefore, the long-time behavior of the Cesaro mean of $u_{\varepsilon}^{E,+}(x,t)$
for the class (C3) is 
\[
M_{t}(u_{\varepsilon}^{E,+}(x,\cdot))\sim(1-\varepsilon)\exp\left(-C\theta_{\varepsilon}^{+}\log(t)^{-1-s}\right)+\varepsilon,\;t\to\infty.
\]
The propagation of the front $x_{\varepsilon,\beta}^{+}(t)$ of $M_{t}(u_{\varepsilon}^{E,+}(x,\cdot))$
of the level $\beta\in(\varepsilon,1)$ is computed solving the following
equation for $x_{\varepsilon,\beta}^{+}(t)$ 
\[
(1-\varepsilon)\exp\left(-C\frac{x_{\varepsilon,\beta}^{+}(t)-x_{\varepsilon}^{+}}{v}\log(t)^{-1-s}\right)+\varepsilon=\beta.
\]
It is easy to find that 
\[
x_{\varepsilon,\beta}^{+}(t)=\log\left(\frac{1-\varepsilon}{\beta-\varepsilon}\right)\frac{v}{C}\log(t)^{1+s}+x_{\varepsilon}^{+}=:C_{+}\log(t)^{1+s}+x_{\varepsilon}^{+}.
\]
Hence, the propagation of the front is $x_{\varepsilon,\beta}^{+}(t)\sim C_{+}\log(t)^{1+s}$
as $t\to\infty$. 

The results for the class (C3) are now stated in the next theorem.
\begin{thm}
\label{thm:FPropagation-C3}Let $u(x,t)=\psi(x-vt)$ be a traveling
wave with constant speed $v$ with two-side estimate, for any $x\in\mathbb{R}$,
$t\ge0$ and $\varepsilon>0$ 
\[
u_{\varepsilon}^{-}(x,t)\le u(x,t)\le u_{\varepsilon}^{+}(x,t).
\]
The subordination $u^{E}(x,t)$ of $u(x,t)$ with the density
$G_{t}(\tau)$ (corresponding to the class (C3)) has Cesaro mean $M_{t}\big(u^{E}(x,\cdot)\big)$
with a propagation of the front $x_{\beta}(t)$ of the level $\beta$
that satisfies the two-side estimate, with $C_{-},C_{+},s>0$ 
\[
C_{-}\log(t)^{1+s}\le x_{\beta}(t)\le C_{+}\log(t)^{1+s}\quad \mathrm{as}\quad t\to\infty.
\]
Then we have $x_{\beta}(t)\sim C\log(t)^{1+s}$ as $t\to\infty$ with $C>0$.
\end{thm}

\subsection*{Acknowledgments}

This work has been partially supported by Center for Research in Mathematics
and Applications (CIMA) related with the Statistics, Stochastic Processes
and Applications (SSPA) group, through the grant UIDB/MAT/04674/2020
of FCT-Funda{\c c\~a}o para a Ci{\^e}ncia e a Tecnologia, Portugal.

The financial support by the Ministry for Science and Education of Ukraine through 
Project 0119U002583 is gratefully acknowledged.


\end{document}